\documentclass[12pt]{amsart}
\usepackage{amssymb,amsmath,amsfonts,latexsym}
\usepackage{bm, enumerate}
\usepackage{hyperref, xcolor}
\usepackage{tikz-cd}

\newcommand{\newsection}[1]{\setcounter{equation}{0} \section{#1}}
\newcommand{\bea}{\begin{eqnarray}}
\newcommand{\eea}{\end{eqnarray}}

\newcommand{\clb}{\mathcal{B}}

\newcommand{\cld}{\mathcal{D}}

\newcommand{\clh}{\mathcal{H}}
\newcommand{\cli}{\mathcal{I}}
\newcommand{\clk}{\mathcal{K}}

\newcommand{\D}{\mathbb{D}}

\newcommand{\C}{\mathbb{C}}

\def\textmatrix#1&#2\\#3&#4\\{\bigl({#1 \atop #3}\ {#2 \atop #4}\bigr)}
\def\dispmatrix#1&#2\\#3&#4\\{\left({#1 \atop #3}\ {#2 \atop #4}\right)}
\newcommand{\be}{\begin{equation}}
\newcommand{\ee}{\end{equation}}
\newcommand{\ben}{\begin{eqnarray*}}
\newcommand{\een}{\end{eqnarray*}}

\newcommand{\NI}{\noindent}

\newcommand{\bi}{\begin{itemize}}
\newcommand{\ei}{\end{itemize}}

\newcommand\la{\langle}
\newcommand\ra{\rangle}
\newtheorem{Theorem}{\sc Theorem}[section]
\newtheorem{Lemma}[Theorem]{\sc Lemma}
\newtheorem{Proposition}[Theorem]{\sc Proposition}
\newtheorem{Corollary}[Theorem]{\sc Corollary}
\newtheorem{Definition}[Theorem]{\sc Definition}
\newtheorem{Question}{\sc Question}
\newtheorem{ass}[Theorem]{\sc Assumption}

\theoremstyle{definition}
\newtheorem{Example}[Theorem]{\sc Example}
\newtheorem{Remark}[Theorem]{\sc Remark}

\newtheorem{Note}[Theorem]{\sc Note}

\newcommand{\bt}{\begin{Theorem}}
\def\beginlem{\begin{Lemma}}
\def\beginprop{\begin{Proposition}}
\def\begincor{\begin{Corollary}}
\def\begindef{\begin{Definition}}
\def\beginexamp{\begin{Example}}
\def\beginrem{\begin{Remark}}
\def\beginq{\begin{Question}}
\def\beginass{\begin{ass}}
\def\beginnote{\begin{Note}}
\newcommand{\et}{\end{Theorem}}
\def\endlem{\end{Lemma}}
\def\endprop{\end{Proposition}}
\def\endcor{\end{Corollary}}
\def\enddef{\end{Definition}}
\def\endexamp{\end{Example}}
\def\endrem{\end{Remark}}
\def\endq{\end{Question}}
\def\endass{\end{ass}}
\def\endnote{\end{Note}}

\numberwithin{equation}{section}

\begin{document}

\title[Zero Products of Toeplitz operators on annulus]{Zero Products of Toeplitz operators on the Hardy and Bergman spaces over an annulus}

\author[Das]{Susmita Das}
\address{Indian Institute of Science, Department of Mathematics, Bangalore, 560012,
India}
\email{susmitadas@iisc.ac.in, susmita.das.puremath@gmail.com}

\author[Narayanan]{E. K. Narayanan}
\address{Indian Institute of Science, Department of Mathematics, Bangalore, 560012,
India}
\email{naru@iisc.ac.in}

\date{\today}
\subjclass[2020]{47B35, 30H10, 30H20}

\keywords{Toeplitz operators, Hardy space, Bergman space, Mellin transform}

\begin{abstract}
We study the zero product problem of Toeplitz operators on the Hardy space and Bergman space over an annulus. Assuming a condition on the Fourier expansion of the symbols, we show that there are no zero divisors in the class of Toeplitz operators on the Hardy space of the annulus. Using the reduction theorem due to Abrahamse, we characterize compact Hankel operators on the Hardy space of the annulus, which also leads to a zero product result. Similar results are proved for the Bergman space over the annulus.
\end{abstract}

\maketitle

\newsection{Introduction}\label{sec: intro}
Let $\mathbb D$ be the open unit disc in $\mathbb C$ and $H^2(\mathbb D)$ be the Hardy space over $\mathbb D.$ For $\varphi \in L^\infty (\mathbb T)$ (where $\mathbb T$ is the unit circle), the Toeplitz operator $T_\varphi$ with symbol $\varphi$ is defined to be $$T_\varphi f = P(\varphi f)$$ where $P$ is the orthogonal projection from $L^2(\mathbb T)$ onto $H^2(\mathbb D).$ The algebraic properties of these operators were studied by Brown and Halmos in their seminal paper \cite{Brown-Halmos}. Among other results, one of the important result they established was that there are no zero divisors for the class of Toeplitz operators. In other words, if $T_\varphi T_\psi = 0,$ then either $\varphi$ or $\psi$ is identically zero. This result has attracted a lot of attention in the past. In particular, there have been attempts to extend this result to other spaces, like the Bergman space over the $\mathbb D,$ and to spaces over domains in higher dimensions. Interestingly, the zero product theorem for the Bergman space in full generality is still open even for the unit disc. In \cite{Ahern-Cuckovic}, Ahern and Cuckovic proved that if the symbols are bounded harmonic functions on the disc, then the zero product theorem is true. Notice that if $\varphi$ is a bounded function on $\mathbb D,$ it admits a polar decomposition $$\varphi(re^{i\theta}) = \sum_{k=-\infty}^{\infty}~\varphi_k(r)~e^{ik\theta},$$ where $\varphi_k(r)$ are the Fourier coefficients of the function $e^{i\theta} \to \varphi(r e^{i\theta}).$ A zero product theorem for Toeplitz operators on the Bergman space over the disc was proved in \cite{LRY} under some assumptions on the polar decomposition of one of the symbols. More precisely, assume that $\psi \in L^\infty(\mathbb D)$ and $\varphi \in L^\infty (\mathbb D),$ with the polar decompositon $$\varphi (re^{i\theta}) = \sum_{k=-\infty}^N~\varphi_k(r)~e^{i k \theta}$$ where $N$ is a positive integer. Assume that $n_0$ is the smallest integer such that $\widehat{\varphi_N}(2n+N+2) \neq 0$ for all $n \geq n_0,$ where $\widehat{\varphi_N}$ is the Mellin transform defined by $$\widehat{\varphi_N}(z) = \int_0^1~\varphi_N(r)~r^{z-1}~dr,$$ then $T_\psi T_\varphi= 0$ implies $\psi = 0.$

\medskip
Our aim in this paper is to prove similar results for the Toeplitz operators defined on the Hardy space and the Bergman space over the annulus $$A = A_{1, R} = \{z \in \mathbb C:~R < |z| < 1 \}.$$ While we follow the methods in \cite{LRY}, we also bring in a powerful theorem, namely the reduction theorem, due to Abrahamse \cite{Abrahamse-Toeplitz} to deal with these questions. The reduction theorem allows us to reduce some of the problems for Toeplitz operators on general multi-connected domains to that of the unit disc. Crucially using this theorem, we also provide a characterization of the compactness of Hankel operators on the annulus, thus establishing an analogue of Hartman's theorem.

\medskip
In the rest of this section, we recall the Hardy space over $A= A_{1, R} = \{z \in \mathbb C:~R< |z| < 1\},$ and some necessary details. This space was introduced and studied in detail by Sarason in \cite{SarasonHp}. We denote by $\partial A$ the boundary of the annulus $A_{1,R}=\{z\in\C: R<|z|<1\}$. Then $\partial A=C\cup C_0$, where $C=\{z\in\C: |z|=1\}$ and $C_0=\{z\in\C: |z|=R\}.$ Let $\text{ Hol}(A)$ be the set of all functions holomorphic on $A_{1,R}$. To define Hardy space $H^2(\partial A)$ on the boundary $\partial A$ of the annulus, as a subspace of $L^2(\partial A)$, we need to introduce the measure, norm, and inner product on $L^2(\partial A)$.

\begin{Definition}\label{H 1}
A subset $E$ of $\partial A$ is called measurable if $\{a\in [0, 2\pi): e^{ia}\in E\}$ and $\{b\in [0, 2\pi): Re^{ib}\in E\}$ are both Borel subsets of $\mathbb{R}$.
\end{Definition}
\NI Let $\sigma$ be the measure defined on $\partial A$ obtained by summing the Lebesgue measure on each component of $\partial A$ and normalised so that $\sigma(\partial A)=2$. More precisely, for $E\subseteq \partial A$ measurable, we define
$$\sigma (E)=\frac{1}{2\pi}\Big(\big(\mu\{a\in[0, 2\pi): e^{ia}\in E\} \big )+\big(\mu\{b\in[0, 2\pi): Re^{ib}\in E\} \big )\Big),$$
where $\mu$ denotes the Lebesgue measure on $\mathbb{R}$.

With this measure $\sigma$, we define the space $L^2(\partial A)$, as the space of all $\sigma$-measurable square integrable functions as follows:
$$L^2(\partial A)=\{f:\partial A\longrightarrow \C: \|f\|_{\partial A}<\infty\}$$
where
$$\|f\|_{\partial A}^2=\frac{1}{2\pi}\int_{0}^{2\pi}|f(e^{it})|^2 dt+\frac{1}{2\pi}\int_{0}^{2\pi}|f(Re^{it})|^2dt,$$
and for $f, g\in L^2(\partial A)$, the corresponding inner product is given by
\[
\begin{split}
\la f,g\ra_{\partial A} & = \int_{\partial A}f\overline{g}~d\sigma \\
&=\frac{1}{2\pi}\int_{0}^{2\pi}f(e^{it})\overline{g(e^{it})}dt+\frac{1}{2\pi}\int_{0}^{2\pi}f(R e^{it})\overline{g(R e^{it})}dt.
\end{split}
\]

The Hardy space $H^2(\partial A)$ is defined to be the closure in $L^2(\partial A)$ of rational functions on $\mathbb C,$ having no poles in $\overline{A}.$
Recall that (see \cite{XH}), the set $\{e_n(z)\}_{n\in\mathbb{Z}}$ forms an orthonormal basis for $H^2(\partial A)$ where
\begin{equation}\label{standard annulus basis}
e_n(z)=\frac{1}{\sqrt{1+R^{2n}}}z^n, \quad z\in\partial A.
\end{equation}
The orthogonal complement of $H^2(\partial A)$ in $L^2(\partial A)$ (see \cite{XH}) is the closed subspace, $\overline{\text{span}}\{ f_n, n\in \mathbb{Z}\},$ where the functions $f_n$ are defined by
\begin{equation}\label{ortho-complement}
f_n(z) =
\begin{cases}
\frac{R^n}{\sqrt{1+R^{2n}}}z^n, & \mbox{if } |z|=1 \\
\frac{-1}{R^n\sqrt{1+R^{2n}}}z^n, & \mbox{if } |z|=R.
\end{cases}
\end{equation}

To study the Toeplitz operator $T_f$ on $H^2(\partial A)$, we need the following definition of Fourier coefficients of $f$ for the outer and inner components $C$ and $C_0$ respectively of $\partial A$.

\begin{Definition}\label{F coefficients}
Let $f\in L^2(\partial A)$. For $n\in \mathbb{Z}$, the $n$-th pair of Fourier coefficients of $f$ denoted by $\widehat{f_C}(n)$ and $\widehat{f_{C_0}}(n)$ respectively and are defined by
$$
\widehat{f_C}(n)=\frac{1}{2\pi}\int_{0}^{2\pi}f(e^{it})e^{-int}dt
$$
$$
\widehat{f_{C_0}}(n)=\frac{1}{2\pi}\int_{0}^{2\pi}f(Re^{it})e^{-int}dt.
$$
\end{Definition}\label{Toeplitz def}

\newsection{Toeplitz operators on $H^2(\partial A)$ }\label{sec: hardy}
We define $L^{\infty}(\partial A)$ to be the space of $\sigma$-measurable essentially bounded functions on $\partial A$.~ Further, $H^{\infty}(\partial A)$ is defined to be the space of all functions in $H^2(\partial A)$ that are also in $L^{\infty}(\partial A)$. Corresponding to $f\in L^{\infty}(\partial A)$, the Toeplitz operator $T_f$ is defined as follows:
\begin{Definition}
Let $P_R$ denote the orthogonal projection of $L^2(\partial A)$ onto $H^2(\partial A)$ and $f\in L^{\infty}(\partial A)$. Then $T_f: H^2(\partial A)\longrightarrow H^2(\partial A)$ is defined by $T_fh=P_R(fh)$, for all $h\in~ H^2(\partial A)$.
\end{Definition}
For $f\in L^{\infty}(\partial A)$ $T_f$ is always bounded (\cite{XH}). A simple computation reveals (\cite{XH},  page 51),
\begin{equation}\label{Toeplitz matrix}
\la T_f e_k, e_j\ra_{\partial A}=\frac{1}{\sqrt{1+R^{2j}}\sqrt{1+R^{2k}}}\big(\widehat{f_C}(j-k)+R^{j+k}\widehat{f_{C_0}}(j-k)\big).
\end{equation}
Equation \eqref{Toeplitz matrix} helps us to write the matrix representation $[T_f]$ of the Toeplitz operator $T_f$ with respect to the orthonormal basis $\{e_n\}_{n\in\mathbb{Z}}$ on $H^2(\partial A)$. Indeed, if $[T_f]=[a_{j,k}]_{j,k=-\infty}^{\infty}$,~ then
\begin{equation*}
[T_f] = \begin{bmatrix}
\vdots & \vdots & \vdots & \vdots & \vdots &\vdots &\vdots
\\
\vdots & \frac{\widehat{f_C}(0)+R^{-4}\widehat{f_{C_0}}(0)}{1+R^{-4}} & a_{-2, -1} & a_{-2,0} & a_{-2,1} & a_{-2,2} & \vdots
\\
\vdots & a_{-1,-2} & \frac{\widehat{f_C}(0)+R^{-2}\widehat{f_{C_0}}(0)}{1+R^{-2}} & a_{-1,0} & a_{-1,1} & a_{-1,2}  & \vdots
\\
\vdots & a_{0,-2} & a_{0,-1} & \frac{\widehat{f_C}(0)+\widehat{f_{C_0}}(0)}{2} & a_{0,1} & a_{0,2} & \vdots
\\
\vdots & a_{1,-2} & a_{1,-1} & a_{1,0} & \frac{\widehat{f_C}(0)+R^{2}\widehat{f_{C_0}}(0)}{1+R^2} & a_{1,2} & \vdots
\\
\vdots  & a_{2,-2} & a_{2,-1} & a_{2,0} & a_{2,1} & \frac{\widehat{f_C}(0)+R^4\widehat{f_{C_0}}(0)}{1+R^4} & \vdots
\\
\vdots & \vdots & \vdots &\vdots &\vdots &\vdots  & \ddots
\end{bmatrix},
\end{equation*}
where
\begin{equation}\label{T entries}
a_{j,k}=\frac{1}{\sqrt{1+R^{2j}}\sqrt{1+R^{2k}}}\big(\widehat{f_C}(j-k)+R^{j+k}\widehat{f_{C_0}}(j-k)\big).
\end{equation}
For $n\in\mathbb{Z}$, we refer to the subdiagonal containing the entries $a_{n,n}=\frac{\widehat{f_C}(0)+R^{2n}\widehat{f_{C_0}}(0)}{1+R^{2n}}$ as the main diagonal of $[T_f]$.
The following lemma will be useful in our context.

\begin{Lemma}\label{T matrix}
$T_f$ is zero if and only if any two columns of $[T_f]$ are zero.
\end{Lemma}
\begin{proof}
It suffices to prove the ``if" part. Let $p\in\mathbb{Z}$, and $C_p$ denote  $p$-th column (that is the column whose entries are $a_{n, p}, n \in \mathbb Z$). For $p \in \mathbb Z,$ let $D_p$ denote the $p$-th subdiagonal (that is the subdiagonal whose entries are $a_{n, n+p}$) of $[T_f]$.

Now, consider two columns $C_r$ and $C_s$ of $[T_f]$ for $r\neq s$. Then, for each $n \in \mathbb Z,$ there exist $m, t\in\mathbb{Z}$ such that
$a_{m,r}\in D_n\cap C_r$ and $a_{t,s}\in D_n\cap C_s$ such that
$$m-r=t-s=n $$
and by \eqref{T entries}, we can write
\begin{equation}\label{amr}
a_{m,r}=\frac{1}{\sqrt{1+R^{2m}}\sqrt{1+R^{2r}}}\big(\widehat{f_C}(n)+R^{m+r}\widehat{f_{C_0}}(n)\big),
\end{equation}
and
\begin{equation}\label{ats}
a_{t,s}=\frac{1}{\sqrt{1+R^{2t}}\sqrt{1+R^{2s}}}\big(\widehat{f_C}(n)+R^{t+s}\widehat{f_{C_0}}(n)\big).
\end{equation}

Since, $r \neq s,$ we have $m+r \neq t+s,$ and since $a_{m,r} = a_{t, s} =0$, we get from \eqref{amr} and \eqref{ats} that
$\widehat{f_C}(n)=\widehat{f_{C_0}}(n)=0$ for every $n$ and so $f = 0.$

\end{proof}

\begin{Remark} Clearly, similar proof works if any two rows are zero.

\end{Remark}

We now prove the following zero product theorem for Toeplitz operators on $H^2(\partial A)$.

\begin{Theorem}\label{Hardy zero product}
Let $g(re^{i\theta})=\sum_{k=-\infty}^{N}g_k (r)e^{ik\theta}$ and $f(re^{i\theta})=\sum_{k=-\infty}^{N'}f_k(r)e^{ik\theta}$, for $r = R, 1$ and $N, N'\in\mathbb{Z},$ be two functions of $L^{\infty}(\partial A)$. Then $T_fT_g=0$ implies $f=0$ or $g=0$.
\end{Theorem}
\begin{proof}
Let $g\neq 0$ and assume without loss of generality, at least one of the Fourier coefficients $\widehat{g_C}(N)$ or $\widehat{g_{C_0}}(N)$ is nonzero.
Then, with respect to $\{e_n\}_{n\in\mathbb{Z}}$, the matrix $[T_g]$ of $T_g$ has an upper triangular form, as the $(j, k)$-th entry $a_{jk}$ is a combination of $\widehat{g_C}(j-k)$ and $\widehat{g_{C_0}}(j-k)$ which is zero provided $j-k > N.$ Notice that the first nonzero subdiagonal from the bottom left corner has entries
$$a_{m,n}=\frac{\widehat{g_C}(N)+R^{m+n}\widehat{g_{C_0}}(N)}{\sqrt{1+R^{2m}}\sqrt{1+R^{2n}}},$$ with $N=m-n$. Moreover, in this subdiagonal, $a_{m,n}$ can vanish at most at one position. Because if there exist distinct $(m_1, n_1)$, $(m_2, n_2)$ such that $a_{m_1, n_1}=a_{m_2, n_2}=0$, where $m_2=m_1+k_1$ and $n_2=n_1+k_1$ for some $k_1(\neq 0)\in\mathbb{Z}$, then
\begin{eqnarray}
\widehat{g_C}(N)+R^{m_1+n_1}\widehat{g_{C_0}}(N) &=& 0 \\
\widehat{g_C}(N)+R^{m_2+n_2}\widehat{g_{C_0}}(N) &=& 0.
\end{eqnarray}

Since $m_2+n_2=m_1+n_1+2k_1\neq m_1+n_1$, it follows that $\widehat{g_C}(N)=\widehat{g_{C_0}}(N)=0$, which contradicts our assumption. Hence we can choose $n_0\in\mathbb{N}$ such that,
\begin{equation}\label{polar Toeplitz}
\widehat{g_C}(N)+R^{2n+N}\widehat{g_{C_0}}(N)\neq 0 \quad \text{ for all } n\geq n_0.
\end{equation}

Now, for any $n\in\mathbb{Z}$ (using the equation \eqref{T entries} for $g$),
\begin{equation}\label{polar Tg}
 T_g(z^n)=\frac{\widehat{g_C}(N)+R^{2n+N}\widehat{g_{C_0}}(N)}{(1+R^{2(n+N)})}z^{n+N}+\sum_{k=-\infty}^{N-1}\frac{\widehat{g_C}(k)+R^{2n+k}\widehat{g_{C_0}}(k)}{(1+R^{2(n+k)})}z^{k+n}.
\end{equation}
Let $T_fT_g=0$. Then for all $n\in \mathbb{Z}$, the equation \eqref{polar Tg} reduces to
\begin{equation}\label{polar Tg zero}
\frac{\widehat{g_C}(N)+R^{2n+N}\widehat{g_{C_0}}(N)}{(1+R^{2(n+N)})}T_f(z^{n+N})+\sum_{k=-\infty}^{N-1}\frac{\widehat{g_C}(k)+R^{2n+k}\widehat{g_{C_0}}(k)}{(1+R^{2(n+k)})}T_f(z^{k+n})=0
\end{equation}
Then  for $n=n_0$, the relation \eqref{polar Toeplitz} and equation \eqref{polar Tg zero} together yield
\begin{equation}\label{n0 plus N case}
T_f(z^{n_0+N})\in\overline{\text{ span}}\{T_f(z^{n_0+N-1}), T_f(z^{n_0+N-2}),\ldots\}.
\end{equation}
Similarly for $n=n_0+1$, it follows by \eqref{polar Toeplitz} and \eqref{polar Tg zero}
\begin{equation}\label{n0 plus N plus one case}
T_f(z^{n_0+N+1})\in\overline{\text{ span}}\{T_f(z^{n_0+N}), T_f(z^{n_0+N-1}), T_f(z^{n_0+N-2}),\ldots\},
\end{equation}
and further \eqref{n0 plus N case} and \eqref{n0 plus N plus one case} together imply
\begin{equation}\label{n0 plus N plus l final}
T_f(z^{n_0+N+1})\in\overline{\text{ span}}\{T_f(z^{n_0+N-1}), T_f(z^{n_0+N-2}),\ldots\}.
\end{equation}
\smallskip
Claim: For $l\geq 0$,
\begin{equation}\label{n0 plus N plus l case}
T_f(z^{n_0+N+l})\in\overline{\text{ span}}\{T_f(z^{n_0+N-1}), T_f(z^{n_0+N-2}),\ldots\}.
\end{equation}
We prove the claim by induction on $l\geq 0$. The proof when $l=0,1$ follows by the equations \eqref{n0 plus N case} and \eqref{n0 plus N plus l final}.
For the induction step, assume the claim to be true for all $0\leq l<m$, for some $m\geq 2$.  Then for $l=m$, it follows by \eqref{polar Toeplitz}, and \eqref{polar Tg zero}
\begin{equation}\label{n0 plus N plus m case}
T_f(z^{n_0+N+m})\in\overline{\text{ span}}\{T_f(z^{n_0+N+m-1}),\ldots, T_f(z^{n_0+N}), T_f(z^{n_0+N-1}),\ldots\},
\end{equation}
and hence the claim follows as by the induction hypothesis,
$$T_f(z^{n_0+N+m-1}),\ldots, T_f(z^{n_0+N})\in \overline{\text{ span}}\{T_f(z^{n_0+N-1}), T_f(z^{n_0+N-2}),\ldots\}.$$

Suppose $f\neq 0$ (equivalently, $T_f\neq 0$). Then
$$\overline{\text{ span}}\{T_f(z^{n_0+N-1}), T_f(z^{n_0+N-2}),\ldots\}\neq 0,$$ otherwise the matrix of $T_f$ will have two columns equal to zero implying $f=0$ (see Lemma \ref{T matrix}). Since $f\neq 0$, there exists an integer $k_0\leq N'$ such that at least one of $\widehat{f_C}(k_0)$ or $\widehat{f_{C_0}}(k_0)$ is nonzero. Note that the matrix of $T_f$ also has an upper triangular form, and the subdiagonal involving  $\widehat{f_C}(k_0), \widehat{f_{C_0}}(k_0)$ can vanish at most at one position.
Since
\begin{equation}\label{Hardy Tf final}
 T_f(z^n)=\sum_{k=-\infty}^{N'}\frac{\widehat{f_C}(k)+R^{2n+k}\widehat{f_{C_0}}(k)}{1+R^{2(n+k)}}z^{k+n},
\end{equation}
we have
\begin{equation}\label{span Tf}
\overline{\text{ span}}\{T_f(z^{n_0+N-1}), T_f(z^{n_0+N-2}),\ldots\}=\overline{\text{ span}}\{z^{n_0+N+N'-1}, z^{n_0+N+N'-2},\ldots\}.
\end{equation}
Now corresponding to $k_0$, there exists $n_{k_0}$ such that
$$n_{k_0}> n_0+N \quad \text{and }\quad n_{k_0}+k_0>n_0+N+N'-1.$$
By \eqref{Hardy Tf final}, we can write
\begin{equation}\label{Hardy Tf last}
T_f(z^{n_{k_0}})=\cdots\cdots+\frac{\widehat{f_C}(k_0)+R^{2n_{k_0}+k_0}\widehat{f_{C_0}}(k_0)}{1+R^{2(n_{k_0}+k_0)}}z^{n_{k_0}+k_0}+\cdots\cdots.
\end{equation}
More generally (again by \eqref{Hardy Tf final}), for all $l'\geq1$
\begin{equation}\label{Hardy Tf l prime}
T_f(z^{n_{k_0}+l'})=\cdots\cdots+\frac{\widehat{f_C}(k_0)+R^{2(n_{k_0}+l')+k_0}\widehat{f_{C_0}}(k_0)}{1+R^{2(n_{k_0}+l')+2k_0}}z^{n_{k_0}+l'+k_0}+\cdots\cdots.
\end{equation}
Clearly, for any $l'\geq 1$
\begin{equation}\label{strict H inq}
n_{k_0}+k_0+l' > n_{k_0}+k_0 > n_0+N+N'-1.
\end{equation}
Now \eqref{n0 plus N plus l case} and \eqref{Hardy Tf final}---\eqref{strict H inq} altogether imply
\begin{eqnarray}
\widehat{f_C}(k_0)+R^{2n_{k_0}+k_0}\widehat{f_{C_0}}(k_0) &=& 0 \\
\widehat{f_C}(k_0)+R^{2(n_{k_0}+l')+k_0}\widehat{f_{C_0}}(k_0) &=& 0,
\end{eqnarray}
which yield $\widehat{f_C}(k_0)=\widehat{f}_{C_0}(k_0)=0$, contradicting our assumption. Hence we must have $f=0$.

For the other case, assume  $T_fT_g=0$ and $f$ is nonzero. If $g\neq0$, then as we have just shown, $f$ must be zero---which is a contradiction. Hence $g$ must be zero.
\end{proof}

Before we go further, let's recall a result from \cite{Abrahamse-Toeplitz} (see Lemma 2.18), which will be used.

\begin{Lemma}\label{abrahamse-lemma}
If $\varphi \in L^\infty (\partial A)$ vanishes on a set of positive measure, but is not identically zero, then Ker~$T_\varphi = 0.$
\end{Lemma}

The above lemma easily implies the following:

\begin{Lemma}\label{lem1}
Let $f,g\in L^{\infty}(\partial A)$ and $T_fT_g=0$. If $fg=0$ on a set $B \subseteq\partial A$ of positive measure, then either $f$ or $g$ is identically zero.
\end{Lemma}
\begin{proof}
If $fg=0$ on $B\subseteq\partial A$ with $\sigma(B)>0$, then there exists $B'\subseteq B$ with $\sigma(B')>0$ such that at least one of $f$ or $g$ vanishes on $B'$. Two cases can arise:

\textbf{Case 1:} $f=0$ on $B'$.
If $f\neq0$ on $\partial A$, then by Lemma \ref{abrahamse-lemma}, $\ker T_f=\{0\}$. Now $T_fT_g= 0$ implies $\text{ Ran} T_g\subseteq \ker T_f=\{0\}$. Hence $T_g=0$ and consequently $g=0$.

\textbf{Case 2:} $g=0$ on $B'$. Since $(T_fT_g)^*=T_{\bar{g}}T_{\bar{f}}$, it follows by case $1$, $\bar{f}=f=0$.
\end{proof}

\subsection{Zero product based on Reduction Theorem}\label{sec: prep res}
Our goal is to obtain some results on the zero product of Toeplitz operators using Abrahamse's reduction theorem. We briefly recall the settings in \cite{Abrahamse-Toeplitz} (see Part III). As earlier, let $A = A_{1, R}$ stand for the annulus $\{z:~R <|z| < 1\}.$ Boundary, $\partial A$ consist of the two circles $C = \{z:~|z| =1\}$ and $C_0 = \{z:~|z| = R\}.$ Interior of $C$ is the unit disc $\mathbb D,$ and let us denote the exterior of $C_0$ including the point $\infty$ by $D_0.$ Thus $D_0 = \{z:~|z| > R\} \cup \{\infty\}.$
Then, by the Caratheodory extension of the Riemann mapping theorem, we get two homeomorphisms $\pi$ and $\pi_0,$ mapping $\mathbb D \cup \mathbb T$ onto $\mathbb D \cup C$ and $D_0 \cup C_0$ respectively, which are conformal equivalences between the interiors. Clearly, we can take $$\pi(z) =z,~\text{and}~\pi_0(z) = R/z.$$
\NI Associated with the function $\phi\in L^{\infty}(\partial A)$ are the functions $\phi_C(z)=\phi \circ \pi(z) = \phi(z)$ and $\phi_{C_0}(z) = \phi \circ {\pi_0}(z) = \phi(R/z),$ in $L^{\infty}(\mathbb T).$ The reduction theorem relates the Toeplitz operator $T_\phi$ with the Toeplitz operators $T_{\phi_C}$ and $T_{\phi_{C_0}}$ on $H^2(\mathbb D)$.

Let $\cli(H^2(\partial A))$ be the $C^*$-algebra generated by $\{T_{\phi}: \phi\in L^{\infty}(\partial A)\},$ and let $\cli(H^2(\mathbb D))$ be the $C^*$-algebra on $H^2(\mathbb D)$ generated by $\{T_{\phi}: \phi\in L^{\infty}(\mathbb T)\}$. For any Hilbert space $\clh$, let $\clb(\clh)$ be the Banach algebra of all bounded operators, $\clk(\clh)$ be the closed ideal of compact operators, and for $T\in\clb(\clh)$, $[T]$ be the coset $T+\clk(\clh)$. Two operators $S$ and $T$ in the same coset are said to be equivalent modulo the compact operators, denoted $S\equiv T$. With these notions, we now state the reduction theorem for Toeplitz operators (\cite{Abrahamse-Toeplitz}, see Theorem 3.1).

\begin{Theorem}\label{Th1}
There is a $*$-isometric isomorphism between the $C^*$-algebras $$\cli(H^2(\partial A))/\clk(H^2(\partial A)) ~\text {and}~ \cli(H^2(\mathbb D))/\clk(H^2(\mathbb D)) \oplus \cli(H^2(\mathbb D))/\clk(H^2(\mathbb D)) $$ which takes $[T_{\phi}]$ to $[T_{\phi_C}] \oplus [T_{\phi_{C_0}}]$.
\end{Theorem}
We shall also need some details about Hankel operators that we recall now. Let the domain $\cld$ stand for $\D$ or $A_{1,R}$ and let $\partial\cld$ be its boundary. The Hankel operator is defined as follows:

\subsection{Hankel operators}
Let $\phi\in L^{\infty}(\partial\cld)$ and $P_R: L^2(\partial\cld)\rightarrow H^2(\partial\cld)$ be the orthogonal projection. Then the Hankel operator $H_{\phi}: H^2(\partial\cld)\rightarrow H^2(\partial\cld)^{\perp}$ is defined by
$$H_{\phi}(f)=(I-P_R)\phi f, \quad\text{ for all }  f\in H^2(\partial\cld).$$

For $\phi \in L^{\infty}(\partial\cld)$, $H_\phi$ is bounded. We need the following lemma from \cite{D Zheng}( see Lemma 1 in \cite{D Zheng}) and \cite{XH} (see the proof of Lemma 3.2.1 in \cite{XH}).

\begin{Lemma}\label{Toeplitz hankel equality}
For $\phi, \psi\in L^{\infty}(\partial\cld)$, $T_{\phi\psi}=T_{\phi} T_{\psi}+H_{\overline{\phi}}^*H_{\psi}$.
\end{Lemma}

For $\phi,\psi\in L^{\infty}(\partial A)$ we have $\phi_{C},~ \phi_{C_0},~\psi_{C}, ~\psi_{C_0} \in L^{\infty}(\mathbb T)$ and by Theorem \ref{Th1} above
$$[T_{\phi}]\longrightarrow [T_{\phi_C}]\oplus[T_{\phi_{C_0}}],\quad [T_{\psi}]\longrightarrow [T_{\psi_C}]\oplus[T_{\psi_{C_0}}] \text{ and}$$
\begin{equation}\label{reduc f}
[T_{\phi}][T_{\psi}]\equiv [T_{\phi}T_{\psi}]\equiv [T_{\phi_C}T_{\psi_C}]\oplus [T_{\phi_{C_0}}T_{\psi_{C_0}}].
\end{equation}
By Lemma \ref{Toeplitz hankel equality}, the relation \eqref{reduc f} further reduces to
\begin{equation}\label{reduc s}
[T_{\phi}T_{\psi}]\equiv [T_{\phi_C \psi_C} - H_{\overline{\phi_C}}^{\ast}H_{\psi_C}] \oplus  [T_{\phi_{C_0}\psi_{C_0}} - H_{\overline{\phi_{C_0}}}^{\ast}H_{\psi_{C_0}}].
\end{equation}
Again for $T_{\phi\psi}$ with $\phi,\psi\in L^{\infty}(\partial A)$, $[T_{\phi\psi}]\equiv [T_{(\phi\psi)_C}]\oplus [T_{(\phi\psi)_{C_0}}]$ (by Theorem \ref{Th1}) and hence
\begin{equation}\label{reduc t}
[T_{\phi\psi}]\equiv [T_{\phi_C\psi_C}]\oplus [T_{\phi_{C_0}\psi_{C_0}}],
\end{equation}
where the second relation follows from $(\phi\psi)_C= \phi_C \psi_C,$ and $(\phi \psi)_{C_0} = \psi_{C_0} \psi_{C_0}.$

Now by \eqref{reduc s} and \eqref{reduc t} $$[T_{\phi}T_{\psi}]\equiv [T_{\phi\psi}] ~\text{if and only if}~ H_{\overline{\phi_C}}^*H_{\psi_C}~\text{and}
H_{\overline{\phi_{C_0}}}^*H_{\psi_{C_0}}~\text{are compact}.$$

\NI Note that, if $H_{\overline{\phi_C}}^*H_{\psi_C}~\text{and }
H_{\overline{\phi_{C_0}}}^*H_{\psi_{C_0}}$ are compact and $T_{\phi}T_{\psi}=0$, then by Lemma \ref{Toeplitz hankel equality}, $T_{\phi\psi}$ is also compact. Since only compact Toeplitz operators are the zero operators on the Hardy space over any domain (\cite{Abrahamse-Toeplitz}, Corollary 2.12), we have $\phi\psi=0$ on $\partial A$ and further by Lemma \ref{lem1}, either $\phi\equiv 0$ or $\psi\equiv 0$. Hence, a zero-product theorem can be obtained via the characterization of compact Hankel operators on the annulus. The following theorem of Hartman (see \cite{Hartman-original} and \cite{Hartman}) on compact Hankel operators for disc plays a crucial role in this sequel.

\begin{Theorem}\label{Hartman}
$H_{\phi}: H^2(\D)\rightarrow H^2(\D)^{\perp}$ is compact if and only if $\phi\in H^{\infty}+C$ on $\mathbb{T}$, where
$$H^{\infty}+C=\{f+g:  f\in H^{\infty}(\mathbb T), g\in C(\mathbb T)\},$$ $C(\mathbb T)$ being the set of all continuous functions on $\mathbb T$.
\end{Theorem}

We now discuss the compactness of the Hankel operator on the annulus $A_{1,R}$. Recall that the boundary $\partial A$ of $A_{1,R}$ consist of the circles $C$ and $C_0$ with $\D$ and $D_0$ being corresponding interior and exterior, including the point at $\infty$ respectively as defined earlier.

Let $A$ and $A_0$ be the algebra of continuous functions on $\mathbb D \cup C$ and $D_0 \cup C_0$ which are holomorphic in $\mathbb D$ and $D_0$ respectively. Also, let $Y$ and $Y_0$ be the closure of $A$ and $A_0$ in $L^2(C)$ and $L^2(C_0)$ respectively. Clearly, $Y= H^2(C)$ and $Y_0=H^2(C_0)$.
Since $L^2(C), L^2(C_0)\subseteq L^2(\partial A)$ are closed, the subspaces $Y$, and $Y_0$ are also closed in $L^2(\partial A)$. Note that, $L^2(\partial A)= H^2(\partial A)\oplus H^2(\partial A)^{\perp}$ and for $\phi\in L^{\infty}(\partial A)$, $H_{\phi}: H^2(\partial A)\rightarrow H^2(\partial A)^{\perp}$ is compact if and only if the operator $ \widetilde{H_{\phi}}=\begin{pmatrix}0 & 0  \\H_{\phi} & 0\end{pmatrix}$ on $H^2(\partial A)\oplus H^2(\partial A)^{\perp}$ is compact. Note by Lemma 3.10 in \cite{Abrahamse-Toeplitz}, $H^2(\partial A)=Y+Y_0$ and hence by Lemma 3.9 in \cite{Abrahamse-Toeplitz}, $\widetilde{H_{\phi}}$ is compact if $\widetilde{H_{\phi}}P_Y$ and $\widetilde{H_{\phi}}P_{Y_0}$ are compact, where $P_Y, P_{Y_0}$ are orthogonal projections of $L^2(\partial A)$ onto $Y$ and $Y_0$ respectively. This is equivalent to $H_{\phi}P_Y$ and $H_{\phi}P_{Y_0}$ being compact.

Again for $\phi\in L^{\infty}(\partial A)$, $\phi\circ\pi, \phi\circ\pi_0\in L^{\infty}(\mathbb{T})$ and one can consider the Hankel operators $H_{\phi}^Y:Y\rightarrow Y^{\perp}(=L^2(C)\ominus Y)$ and $H_{\phi}^{Y_0}:Y_0\rightarrow Y_0^{\perp}(=L^2(C_0)\ominus Y_0)$ defined by $H_{\phi}^Y(f)=P_{Y^{\perp}}\phi f$, $f\in Y$ and $H_{\phi}^{Y_0}(f)=P_{Y_0^{\perp}}\phi g$ $g\in Y_0$; $P_{Y^{\perp}}, P_{Y_0^{\perp}}$ being the orthogonal projections of $L^2(C), L^2(C_0)$ onto $Y^{\perp}$ and $Y_0^{\perp}$ respectively. Our goal is to relate the components $H_{\phi}P_Y, H_{\phi}P_{Y_0}$ with $H_{\phi}^Y$ and $H_{\phi}^{Y_0}$ respectively to derive a compactness criterion for $H_{\phi}$. In what follows, we will deal with $H_{\phi}P_{Y_0}, H_{\phi}^{Y_0}$ only, as $H_{\phi}P_{Y}, H_{\phi}^Y$ can be shown related exactly in the same way.

Corresponding to the homeomorphisms $\pi,\pi_0$, define the maps $\widetilde{\pi}$ and $\widetilde{\pi_0}$ from $L^2(C)$ and $L^2(C_0)$ respectively to $L^2(\mathbb T)$ by $\widetilde{\pi}(f) = f \circ \pi$ and $\widetilde{\pi_0}(g) = g \circ \pi_0.$ Then $\{\widetilde{e_n} \circ \pi^{-1} \}_{n \in \mathbb Z}$ and $\{\widetilde{e_n} \circ \pi_0^{-1}\}_{n \in \mathbb Z}$ form orthonormal bases for $L^2(C)$ and $L^2(C_0)$ respectively, where $\{\widetilde{e_n}\}_{n \in \mathbb Z},$ with $\widetilde{e_n}(z) = z^n,~z \in \mathbb T,$ is the standard orthonormal basis for $L^2(\mathbb T)$. Note that, the sets $\{\widetilde{e_n} \circ \pi^{-1} \}_{n \in \mathbb Z_{\geq 0}}$ and
$\{\widetilde{e_n} \circ \pi_0^{-1} \}_{n \in \mathbb Z_{\geq 0}}$ form orthonormal bases for $Y$ and $Y_0$ respectively.

Let us define $U: Y^{\perp}\longrightarrow H^2(\mathbb D)^{\perp}~\text{and}~U_0: Y_0^{\perp} \longrightarrow H^2(\mathbb D)^{\perp}$ by
\begin{equation}\label{define U}
U \big ( \sum_{n \leq -1}\langle f, ~\widetilde{e_n} \circ \pi^{-1} \rangle \widetilde{e_n} \circ \pi^{-1} \big ) = \sum_{n \leq -1} \langle f \circ \pi, \widetilde{e_n} \rangle \widetilde{e_n}, \quad f\in Y
\end{equation}
\begin{equation}\label{define U knot}
U_0 \big ( \sum_{n \leq -1}\langle f, ~\widetilde{e_n} \circ \pi_0^{-1} \rangle \widetilde{e_n} \circ \pi_0^{-1} \big ) = \sum_{n \leq -1} \langle f \circ \pi_0, \widetilde{e_n} \rangle \widetilde{e_n}, \quad f\in Y_0.
\end{equation}

Clearly, $U$ and $U_0$ are unitary maps. We show that the following diagram is commutative.
\[\begin{tikzcd}\label{dia}
Y_0 \arrow{r}{\widetilde{\pi_0}} \arrow[swap]{d}{H_{\phi}^{Y_0}} & H^2(\D) \arrow{d}{ H_{\phi_{C_0}}(=H_{\phi\circ\pi_0})} \\
Y_0^{\perp} \arrow{r}{U_0} & H^2(\D)^{\perp}
\end{tikzcd}
\]
Indeed for $f\in Y_0$, $\widetilde{\pi_0}(f)=f\circ\pi_0=\tilde{f}$ (say), and hence
\begin{equation}\label{hankel one}
H_{\phi}^{Y_0}\widetilde{\pi_0}^{-1}(\tilde{f})=H_{\phi}^{Y_0}(f)=P_{{Y_0}^{\perp}}(\phi f).
\end{equation}
Since $\{\widetilde{e_n}\circ\pi_0^{-1}\}_{n\leq -1}$ is an orthonormal basis of $Y_0^{\perp}$,
\begin{equation}\label{hankel two}
P_{{Y_0}^{\perp}}(\phi f)=\sum_{n\leq -1}\langle \phi f, \widetilde{e_n}\circ\pi_0^{-1}\rangle \widetilde{e_n}\circ\pi_0^{-1}.
\end{equation}
By equations \eqref{define U knot}, \eqref{hankel one} and \eqref{hankel two},
\begin{equation}\label{U knot and hankel}
U_0H_{\phi}^{Y_0}\widetilde{\pi_0}^{-1}(\tilde{f})=\sum_{n\leq -1}\langle (\phi f)\circ\pi_0, \widetilde{e_n}\rangle \widetilde{e_n}.
\end{equation}
Since $\phi f\circ\pi_0=(\phi\circ\pi_0)(f\circ\pi_0)$ and $\phi\circ\pi_0= \phi_{C_0}$, \eqref{U knot and hankel} becomes
\begin{equation}\label{u knot hankel final}
U_0H_{\phi}^{Y_0}\widetilde{\pi_0}^{-1}(\tilde{f})=\sum_{n\leq -1}\langle \phi_{C_0}\tilde{f}, \widetilde{e_n}\rangle \widetilde{e_n} =H_{\phi_{C_0}}(\tilde{f}).
\end{equation}
\NI Since $U_0$ and $\widetilde{\pi_0}^{-1}$ are invertible, by \eqref{u knot hankel final} $H_{\phi}^{Y_0}$ is compact if and only if $H_{\phi_{C_0}}$ is compact---which is, by the result of Hartman (Theorem \ref{Hartman}) equivalent to $\phi_{C_0}\in H^{\infty}+C$.

Note that, $\phi\in L^{\infty}(\partial A)$ implies $\phi\in L^{\infty}(C_0)$ and Recall:
$$L^2(C_0)=Y_0\oplus Y_0^{\perp}= H^2(C_0)\oplus H^2(C_0)^{\perp}\subseteq H^2(\partial A)\oplus H^2(\partial A)^{\perp}.$$
Then for $f\in Y_0$ and $\phi\in L^{\infty}(\partial A)$, $\phi f\in L^2(C_0)$ and one can write
\begin{equation}\label{h comp 1}
\phi f= y_1\oplus y_2\quad \text{ for some } y_1\in Y_0,\quad y_2\in Y_0^{\perp}.
\end{equation}
Again, as $Y_0^{\perp}\subseteq H^2(\partial A)\oplus H^2(\partial A)^{\perp}$, there exist $x_1\in H^2(\partial A)$ and $x_2\in H^2(\partial A)^{\perp}$ such that
\begin{equation}\label{y2 as direct sum}
H_{\phi}^{Y_0}f=y_2=x_1\oplus x_2.
\end{equation}
On the other hand, for the same $f$ in $Y_0$, it follows by equation \eqref{h comp 1}
\begin{equation}\label{H phi proj on Y knot 1}
H_{\phi}P_{Y_0}(f)=H_{\phi}(f)=P_{H^2(\partial A)^{\perp}}(\phi f)=P_{H^2(\partial A)^{\perp}}(y_1\oplus y_2).
\end{equation}
Since $y_1\in Y_0\subseteq H^2(\partial A)$, equations \eqref{y2 as direct sum} and \eqref{H phi proj on Y knot 1} together imply
\begin{equation}\label{H phi proj on Y knot 2}
H_{\phi}P_{Y_0}(f)=x_2,
\end{equation}
and finally by \eqref{y2 as direct sum} and \eqref{H phi proj on Y knot 2}
\begin{equation}\label{PH=HP}
P_{H^2(\partial A)^{\perp}}H_{\phi}^{Y_0}=H_{\phi}P_{Y_0}.
\end{equation}

Note that, $H_{\phi}^{Y_0}$ can be decomposed as
\begin{equation}\label{H phi o}
H_{\phi}^{Y_0}=P_{H^2(\partial A)}H_{\phi}^{Y_0}\oplus P_{H^2(\partial A)^{\perp}}H_{\phi}^{Y_0}.
\end{equation}
We show that the component $P_{H^2(\partial A)}H_{\phi}^{Y_0}$ can be written as the product $TP_{H^2(\partial A)^{\perp}}H_{\phi}^{Y_0}$, for some suitable $T: H^2(\partial A)^{\perp}\rightarrow H^2(\partial A)$. Then we can conclude through \eqref{PH=HP} and \eqref{H phi o}, $H_{\phi}P_{Y_0}$ is compact if and only if $H_{\phi}^{Y_0}$ is compact.

Since $y_2\in Y_0^{\perp}=L^2(C_0)\ominus H^2(C_0)$, we have $y_2\circ\pi_0\in H^2(\D)^{\perp}$ (see the diagram above) and hence $\overline{y_2\circ\pi_0}\in H^2(\D)$. Since $\pi_0$ is a homeomorphism, $\overline{y_2}\in Y_0\subseteq H^2(\partial A)$ and one can write
\begin{equation}\label{y complex conj}
\overline{y_2}=\sum_{n\in\mathbb{Z}}\langle \overline{y_2}, e_n\rangle_{L^2(\partial A)}e_n,
\end{equation} where $\{e_n\}_{n\in \mathbb Z}$ is as in \eqref{standard annulus basis}.
Then by \eqref{y2 as direct sum} and \eqref{y complex conj}
\begin{equation}\label{y2 annulus}
y_2=\sum_{n\in\mathbb{Z}}\langle e_n,\overline{y_2}\rangle_{L^2(\partial A)}\overline{e_n}=H_{\phi}^{Y_0}f.
\end{equation}
\NI At this point, we need the following lemma:

\begin{Lemma}\label{en bar annulus}
$\overline{e_n}=\frac{2R^n}{1+R^{2n}}e_{-n}+\frac{1-R^{2n}}{1+R^{2n}}f_{-n}$ on $L^2(\partial A)$.
\end{Lemma}
\begin{proof}
Recall that, the set $\{e_n,f_n\}_{n\in\mathbb{Z}}$ is an orthonormal basis of $L^2(\partial A)$ (see \eqref{standard annulus basis} and \eqref{ortho-complement}) and hence for any $n\in \mathbb{Z}$,
\begin{equation}\label{conjugate basis}
\overline{e_n}=\sum_{m\in\mathbb{Z}}\langle\overline{e_n},e_m\rangle_{\partial A}~e_m + \sum_{m\in\mathbb{Z}}\langle\overline{e_n},f_m\rangle_{\partial A}~f_m.
\end{equation}
The proof now follows by a direct computation after substituting  $e_m$ and $f_m$ in \eqref{conjugate basis}.
\end{proof}
Clearly by Lemma \ref{en bar annulus} , the equation \eqref{y2 annulus} further reduces to
\begin{equation}\label{y2 as two pfojc}
y_2=\sum_{n\in\mathbb{Z}}\langle e_n,\overline{y_2}\rangle_{L^2(\partial A)}\Big[\frac{2R^n}{1+R^{2n}}e_{-n}+\frac{1-R^{2n}}{1+R^{2n}}f_{-n}\Big].
\end{equation}
Again, by \eqref{y2 annulus} and \eqref{y2 as two pfojc}
\begin{eqnarray}
P_{H^2(\partial A)}H_{\phi}^{Y_0} f&=& P_{H^2(\partial A)}y_2\ =\ \sum_{n\in\mathbb{Z}}\langle e_n,\overline{y_2}\rangle_{L^2(\partial A)}\Big(\frac{2R^n}{1+R^{2n}}\Big) e_{-n},\ \ \text{and}\label{Py2}\\
P_{H^2(\partial A)^{\perp}}H_{\phi}^{Y_0} f &=& P_{H^2(\partial A )^{\perp}}y_2\ =\ \sum_{n\in\mathbb{Z}}\langle e_n,\overline{y_2}\rangle_{L^2(\partial A)}\Big(\frac{1-R^{2n}}{1+R^{2n}}\Big)f_{-n}.\label{Py2 perp}
\end{eqnarray}
Let us now define $T: H^2(\partial A)^{\perp}\longrightarrow H^2(\partial A)$ by
\begin{equation}\label{define T}
f_{-n}\longrightarrow \frac{2R^n}{1-R^{2n}}e_{-n},\quad \text{ for all } n\in\mathbb{Z}.
\end{equation}
Then we have $T=T_2T_1$, where $T_1: H^2(\partial A)^{\perp}\rightarrow H^2(\partial A)^{\perp}$ and $T_2: H^2(\partial A)^{\perp}\rightarrow H^2(\partial A)$ are defined by
$$T_1(f_n)=\frac{2R^n}{1-R^{2n}}f_n\quad \text{ for all } n\in\mathbb{Z},$$
and $T_2(f_n)=e_n$ $\forall n\in\mathbb{Z}$. Then by \eqref{Py2}, \eqref{Py2 perp}, and \eqref{define T}
\begin{equation}\label{PH phi}
P_{H^2(\partial A)}H_{\phi}^{Y_0} f=P_{H^2(\partial A)} y_2= TP_{H^2(\partial A)^{\perp}} y_2=TP_{H^2(\partial A)^{\perp}}H_{\phi}^{Y_0} f.
\end{equation}
Since $f\in Y_0$ is arbitrary, we have
\begin{equation}\label{PH phi perp}
P_{H^2(\partial A)}H_{\phi}^{Y_0}= TP_{H^2(\partial A)^{\perp}}H_{\phi}^{Y_0}.
\end{equation}
Note that $T_1$ and $T$ are compact. As we mentioned earlier, one can proceed similarly to obtain the same relations and conclusions about $H_{\phi}^Y, H_{\phi}P_Y$. So far discussion summarizes to:

\begin{Theorem}\label{Hankel Annulus}
The operator $H_{\phi}: H^2(\partial A)\longrightarrow H^2(\partial A)^{\perp}$ with $\phi\in L^{\infty}(\partial A)$ is compact if and only if the functions $\phi_C,\phi_{C_0}\in H^{\infty}+C$.
\end{Theorem}
We now state the following zero product theorem for the Toeplitz operators on $H^2(\partial A)$---the proof of which follows by the reduction Theorem \ref{Th1} , Theorem \ref{Hankel Annulus} , Lemma \ref{lem1} , and Lemma \ref{Toeplitz hankel equality} .
\begin{Theorem}\label{0 product thm}
Let $\phi,\psi\in L^{\infty}(\partial A)$ such that $\overline{\phi}_C$ and $\overline{\phi}_{C_0}$(\text{or} $\psi_C$ and $ \psi_{C_0}$) belong to $H^{\infty}+C$. Then  $T_{\phi}T_{\psi}=0$ on $H^2(\partial A),$ implies $\phi=0$ or $\psi=0$.
\end{Theorem}

\newsection{Toeplitz operators on the Bergman space of the annulus}\label{sec: bergman}

In this section we look at the zero product theorem on the Bergman space over an annulus.
Let $R>0$. As earlier, we denote by $A_{1,R}$ the annulus
$$A_{1,R}=\{z\in \C: R<|z|<1 \}.$$ Recall that, the Bergman space $B^2(A_{1,R})$ is the space of all square integrable holomorphic functions on $A_{1,R}$ i.e.,
$$B^2(A_{1,R})=\{f:A_{1,R}\rightarrow \C,\text{ holomorphic and } \int_{A_{1,R}}|f(z)|^2~ dA(z)<\infty\},$$ where $dA(z) = dxdy$ is the area measure. For $f, g\in B^2(A_{1, R})$, the norm and inner product of the space are given by
$$\|f\|_{B^2(A_{1, R})}^2=\frac{1}{2\pi}\int_{A_{1, R}}|f(z)|^2 dA(z),$$
$$\la f,g\ra_{B^2(A_{1, R})}=\frac{1}{2\pi}\int_{A_{1, R}}f(g)\overline{g(z)} dA(z).$$
It is well-known that $B^2(A_{1, R})$ is a closed subspace of $L^2(A_{1, R}, dA)$. Following lemma provides an orthonormal basis for $B^2(A_{1,R})$.

\begin{Lemma}\label{Bergman basis}
The set $\{\sqrt{\frac{2(n+1)}{1-R^{2(n+1)}}}z^n\}_{n\in\mathbb{Z}\setminus\{-1\}}\bigcup \{\frac{z^{-1}}{(\log \frac{1}{R})^{1/2}}\}$ is an orthonormal basis of $B^2(A_{1,R})$.
\end{Lemma}
\begin{proof}
Follows by an easy and straightforward computation.
\end{proof}
We now introduce the Mellin transform, radial, and quasi-homogeneous functions that will be useful in our context.
\begin{Definition}\label{Mellin transform}
The Mellin transform $\widehat{\phi}$ of of a function $\phi\in L^1([R, 1], rdr)$ is defined by
$$\widehat{\phi}(z)=\int_{R}^{1}\phi(r) r^{z-1}dr$$
\end{Definition}

In fact, the Mellin transform is defined for suitable functions defined on $(0, \infty).$ In the above, the function is considered to be zero on $(0, R) \cup (1, \infty).$ For a function $\psi\in L^1([0,1], rdr)$ \big(considered to be zero on $(1,\infty)$\big), it is well-known (\cite{GV}.\cite{LRY},\cite{LSZ}) that the Mellin transform $\widehat{\psi}$ is well-defined on $\{z\in\C: \text{ Re }z\geq 2\}$ and analytic on $\{z\in\C: \text{ Re }z> 2\}$. Also, a function can be determined by the value of a certain number of its Mellin coefficients. The following lemma proved in (\cite{GV}) establishes this:

\begin{Lemma}\label{Mellin-classical result 1}
Let $\phi\in L^1([0,1], rdr)$. If there exist $n_0, p\in\mathbb{Z}$ such that
$$\widehat{\phi}(n_0+pk)=0 \quad \text{ for all } k\in\mathbb{N},$$ then $\phi=0$.
\end{Lemma}

\begin{Definition}\label{radial}
A function $\phi\in L^1(A_{1, R}, dA)$ is called radial if
$$\phi(z)=\phi_1(|z|)\quad R\leq |z|\leq 1,$$ for a function $\phi_1$ on $[R, 1].$
A function $f$ defined on $A_{1, R}$ is said to be quasi-homogeneous of degree $p \in \mathbb{Z}$ if we can write it as $e^{ip\theta}\phi$, where $\phi$ is a radial function.
\end{Definition}
Clearly, a radial function is a quasi-homogeneous function of degree zero. Note that, any function $f\in L^2(A_{1, R})$ has the polar decomposition
$$f(r e^{i\theta})=\sum_{k\in\mathbb{Z}}f_k(r)e^{ik\theta},$$
where $f_k$ are radial in $L^2([R,1], rdr)$. Below, we define the Toeplitz operator on $B^2(A_{1, R})$.
\begin{Definition}
Let $L^\infty(A_{1, R})$ denote the algebra of all essentially bounded functions on $A_{1, R}$ and $P_{B^2(A_{1, R})}$ be the orthogonal projection from $L^2(A_{1, R})$ onto $B^2(A_{1,R})$. Corresponding to $f \in L^\infty(A_{1, R})$, the Toeplitz operator $T_f:B^2(A_{1, R})\rightarrow B^2(A_{1, R})$ is defined by $T_f \varphi =~P_{B^2(A_{1, R})}(f \varphi)$.
\end{Definition}

 A Toeplitz operator $T_f$ where $f$ is quasi-homogeneous, is called a quasi-homogeneous Toeplitz operator. We now determine the Toeplitz operator $T_f$ where $f$ is a quasi-homogeneous symbol given by $f(r e^{i\theta})=f_1(r)e^{ip\theta}, p\in\mathbb{Z}$.  Let us set
\begin{equation}\label{t as basis norm}
t_n=
\begin{cases}
&\frac{1}{(\log\frac{1}{R})^{1/2}} \mbox{ if } n= -1 \\

& \sqrt{\frac{2(n+1)}{1-R^{2(n+1)}}} \mbox{ if } n\neq -1.
\end{cases}
\end{equation}
 Then we have the following lemma:
\begin{Lemma}\label{T operator}
For $n\in\mathbb{Z}$,
$$T_f(z^n)=t_{p+n}^2\widehat{f_1}(p+2n+2)z^{p+n}.$$
\end{Lemma}
\begin{proof}
Note by Lemma \ref{Bergman basis}, and equation (\ref{t as basis norm}), an orthonormal basis of $B^2(A_{1, R})$ is given by $\{t_n z^n\}_{n\in\mathbb{Z}}$. Then for $f=e^{i p\theta}f_1(r)$ with $p\in\mathbb{Z}$, it follows that $\forall n \in \mathbb{Z}$,
\begin{equation}\label{quasi-toeplitz}
T_f(z^n)=P_{B^2(A_{1, R})}(fz^n)=\sum_{m\in\mathbb{Z}}\la f z^n, t_m z^m\ra_{A_{1, R}} t_m z^m, \quad \text{where}
\end{equation}
\[
\begin{split}
\la f z^n, t_mz^m\ra_{A_{1, R}} &=\frac{1}{2\pi}\int_{A_{1, R}} t_m f(z) z^n\overline{z}^m dA(z) \\
&=\frac{1}{2\pi}t_m\int_{R}^{1}f_1(r)r^{m+n+1}dr\int_{0}^{2\pi}e^{i(p+n-m)}d\theta \\
&=t_m\widehat{f_1}(n+m+2)\frac{1}{2\pi}\int_{0}^{2\pi}e^{i(p+n-m)\theta}d\theta, \quad\text{and hence}\\
\end{split}
\]
\[
\la f z^n, t_mz^m\ra_{A_{1, R}}=
\begin{cases}
t_{p+n}\widehat{f_1}(p+2n+2), & \mbox{if } m=p+n \\
0, & \mbox{if }  m\neq p+n.
\end{cases}
\]
Therefore, it follows by \eqref{quasi-toeplitz}
$$T_f(z^n)=t_{p+n}^2\widehat{f_1}(p+2n+2)z^{p+n}.$$
\end{proof}

\begin{Note}
$T_f$ is a forward (backward) shift if $p>0\quad ( p<0)$, and a diagonal operator if $p=0$, i.e., $f$ is radial.
\end{Note}
We now prove the following zero product theorem.
\begin{Theorem}\label{Bergman zero product}
Let $f,g\in L^{\infty} (A_{1,R})$ such that $f(re^{i\theta})=\sum_{k=-\infty}^{M}f_k(r)e^{ik\theta}$ and $g(re^{i\theta})=\sum_{k=-\infty}^{N}g_k(r)e^{ik\theta}$ for some $M,N\in\mathbb{Z}$. Assume $n_0\in\mathbb{Z}$ to be the smallest integer such that $\widehat{g_N}(2n+N+2)\neq 0$ for all $n\geq n_0$. If $T_fT_g=0$ then $f=0$.
\end{Theorem}
\begin{proof}
By Lemma \ref{T operator}, for all $n\in\mathbb{Z}$
\begin{equation}\label{Tg zn}
T_g(z^n)= t_{N+n}^2 \widehat{g_N}(N+2n+2)z^{n+N}+\sum_{k=-\infty}^{N-1}t_{k+n}^2 \widehat{g_k}(k+2n+2)z^{k+n}.
\end{equation}

Since $\widehat{g_N}(2n_0+N+2)\neq 0$,
\begin{equation}\label{berg1}
z^{n_0+N}\in\overline{\text{ span}}\{T_g(z^{n_0}), z^{n_0+N-1}, z^{n_0+N-2},\ldots\}.
\end{equation}
Similarly for $n=(n_0+1)$, $\widehat{g_N}(2n_0+N+4)\neq 0$, and equation \eqref{Tg zn} implies
\begin{equation}\label{berg2}
z^{n_0+N+1}\in\overline{\text{ span}}\{T_g(z^{n_0+1}), z^{n_0+N}, z^{n_0+N-1}, z^{n_0+N-2},\ldots\}.
\end{equation}
Hence \eqref{berg1} and \eqref{berg2} together yield
\begin{equation}\label{berg3}
z^{n_0+N+1}\in\overline{\text{ span}}\{T_g(z^{n_0+1}), T_g(z^{n_0}), z^{n_0+N-1}, z^{n_0+N-2},\ldots\}.
\end{equation}
Proceeding exaclty in the same way, it follows by induction that, for all $l\geq 0$
\begin{equation}\label{Tg zn+N+l}
z^{n_0+N+l}\in\overline{\text{ span}}\{T_g(z^{n_0+l}),\ldots, T_g(z^{n_0}), z^{n_0+N-1}, z^{n_0+N-2},\ldots\}.
\end{equation}
Let $T_fT_g=0$. Then the relation \ref{Tg zn+N+l} reduces to
\begin{equation}\label{Tf span}
T_f(z^{n_0+N+l})\in\overline{\text{ span}}\{T_f(z^{n_0+N-1}), T_f(z^{n_0+N-2}),\ldots\} \text{ for all } l\geq 0.
\end{equation}

We now consider $T_f$. for all $n\in\mathbb{Z}$, Lemma \ref{T operator} implies
\begin{equation}\label{Tf zn}
T_f(z^n)=\sum_{k=-\infty}^{M}t_{k+n}^2\widehat{f_k}(k+2n+2) z^{k+n},
\end{equation}
and hence for all $n\in\mathbb{Z}$,
\begin{equation}\label{Tf M+n}
T_f(z^n)\in\overline{\text{ span}}\{z^{M+n}, z^{M+n-1},\ldots\ldots\}.
\end{equation}
Now for  $l\geq 0$, equations \ref{Tf span} and \ref{Tf M+n} together imply
\begin{equation}\label{Tf=M+n}
T_f(z^{n_0+N+l})\in\overline{\text{ span}}\{z^{M+(n_0+N-1)}, z^{M+(n_0+N-2)}, z^{M+(n_0+N-3)},\ldots\}.
\end{equation}

For $M\in\mathbb{Z}$, there exists $l_M\in\mathbb{N}$ such that
$$M+2(n_0+N+l_M+1)\in\mathbb{N}$$ and
$$M+n_0+N+l_M> M+n_0+N-1$$
Then by \ref{Tf zn},
\begin{equation}\label{T_f lm}
T_f(z^{n_0+N+l_M})=\sum_{k=-\infty}^{M}t_{k+(n_0+N+l_M)}^2\widehat{f_k}\big(k+2(n_0+N+l_M)+2\big) z^{k+(n_0+N+l_M)}.
\end{equation}
Now for any $l\geq 0$, we have
\begin{equation}\label{lm and l}
M+n_0+N+l_M+l>M+n_0+N-1
\end{equation}
and again by \ref{Tf zn}, for all $l\geq 0$,
\begin{equation}\label{T_f final}
T_f(z^{n_0+N+l_M+l})=\sum_{k=-\infty}^{M}t_{k+(n_0+N+l_M+l)}^2\widehat{f_k}(k+2(n_0+N+l_M+l)+2) z^{k+(n_0+N+l_M+l)}.
\end{equation}
Hence it follows by \eqref{Tf=M+n}, \eqref{lm and l}, and \eqref{T_f final}
\begin{equation}\label{fM radial}
\widehat{f_M}(M+2(n_0+N+l_M+l)+2)=0, \quad\forall l\geq 0.
\end{equation}
Now for the radial function $f_M (\in (L^1[R, 1], rdr))$,  $M+2(n_0+N+l_M+1), 2\in\mathbb{N}$ such that
$$\widehat{f_M}(M+2(n_0+N+l_M+1)+2l)=0 \quad\forall l\geq 0 \quad (\text{ by }\eqref{fM radial}).$$
Hence by Lemma \ref{Mellin-classical result 1}, $f_M=0$.
Similarly, for any $k\in(-\infty, M)$ and for any radial function $f_k(r)$, there exists $l_k\in\mathbb{N}$ such that
\[
k+2(n_0+N+l_k+1)\in\mathbb{N}
\]
and
\[
k+n_0+N+l_k>M+n_0+N-1,
\]
and again, by a similar argument as above, we have
\[
\widehat{f_k}(k+2(n_0+N+l_k+1)+2l)=0 \text{ for all } l\geq 0,
\]
and hence by Lemma \ref{Mellin-classical result 1}, $f_k=0$.
Since $k\in (-\infty, M)$ is an arbitrary integer, it follows that $f_k=0$ for all $k\in (-\infty, M)\cap\mathbb{Z}$ and hence $f=0$.
\end{proof}
\vspace{0.1 in}

\textbf{Acknowledgement:}
The research of the first-named author is supported by NBHM Postdoctoral Fellowship at the Indian Institute of Science, Bangalore, India.

\bibliographystyle{amsplain}

\end{document}